\documentclass[12pt]{amsart}
\usepackage[colorlinks=true,citecolor=green,linkcolor=red]{hyperref}
\usepackage{amsmath}
\usepackage{verbatim}
\usepackage{amsfonts}
\usepackage{amssymb}
\usepackage{color}
\usepackage[all]{xy}  
\usepackage{enumerate}
\usepackage[top=1in, bottom=1in, left=0.9in, right=0.9in]{geometry}
\usepackage{mathrsfs}

\usepackage{stmaryrd}
\usepackage{bbold}

\usepackage{cleveref}
\usepackage{soul}

\theoremstyle{definition}
\newtheorem{thm}{Theorem}[section]

\numberwithin{equation}{section}

\theoremstyle{definition}

\newtheorem{example}[thm]{Example}

\newtheoremstyle{TheoremNum}
        {8pt}{8pt}              
        {\upshape}                      
        {}                              
        {\bfseries}                     
        {.}                             
        {.5em}                             
        {\thmname{#1}\thmnote{ \bfseries #3}}
  \theoremstyle{TheoremNum}

\newcommand{\Ext}{\operatorname{Ext}}

\newcommand{\Ass}{\operatorname{Ass}}

\newcommand{\Min}{\operatorname{Min}}

\title{Calculations involving symbolic powers}

\author[Drabkin]{Ben Drabkin}
\address{Department of Mathematics, University of Nebraska, Lincoln, NE 68588-0130, USA}
\email{benjamin.drabkin@huskers.unl.edu}

\author[Grifo]{Elo\'isa Grifo}
\address{Department of Mathematics, University of Virginia, Charlottesville, VA 22904-4135, USA}
\email{eloisa.grifo@virginia.edu}

\author[Seceleanu]{Alexandra Seceleanu}
\address{Department of Mathematics, University of Nebraska, Lincoln, NE 68588-0130, USA}
\email{aseceleanu@unl.edu}

\author[Stone]{Branden Stone}
\address{Dept of Math and Computer Science, Adelphi University, Garden City, NY 11530-0701, USA}
\email{bstone@adelphi.edu}


\begin{document}

\begin{abstract}
	Symbolic powers are a classical commutative algebra topic that relates to primary decomposition, consisting, in some circumstances, of the functions that vanish up to a certain order on a given variety. However, these are notoriously difficult to compute, and there are seemingly simple questions related to symbolic powers that remain open even over polynomial rings. In this paper, we describe a \href{https://faculty.math.illinois.edu/Macaulay2/}{\emph{Macaulay2}} software package that allows for computations of symbolic powers of ideals and which can be used to study the equality and containment problems, among others.
\end{abstract}

\maketitle

\section{Introduction}


Given an ideal $I$ in a Noetherian domain $R$, the \emph{$n$-th symbolic power of I} is the ideal defined by
$$I^{(n)} = \bigcap_{P \in \Ass(I)} \left( I^n R_P \cap R \right).$$
When $I$ itself has no embedded primes, the minimal primes of $I^n$ coincide with the associated primes of $I$, and $I^{(n)}$ as above corresponds to the intersection of the primary components corresponding to minimal primes of $I^n$. In particular, under these circumstances the definition is unchanged if instead we have $P$ ranging over $\Min(I)$. However, if we consider any ideal $I$, with no assumptions on its associated primes, there are two possible definitions of symbolic powers: the one above, and the one given by
$$I^{(n)} = \bigcap_{P \in \Min(I)} \left( I^n R_P \cap R \right).$$
The \texttt{SymbolicPowers} package, in
\href{https://github.com/eloisagrifo/SymbolicPowers}{https://github.com/eloisagrifo/SymbolicPowers}, allows the user to compute the symbolic powers of any ideal over a polynomial ring, taking the first definition as the standard, but allowing the user to take the second definition instead via the option \texttt{UseMinimalPrimes}, which can be used in any method.

Symbolic powers are a classical topic that relates to many subjects within commutative algebra and algebraic geometry, and an active area of current research. If $P$ is a prime ideal in a polynomial ring, the classical Zariski--Nagata \cite{Zariski,Nagata} theorem says that the symbolic powers of $P$ consist of the functions that vanish up to order $n$ in the corresponding variety. Over a perfect field, these coincide with differential powers. For a survey on symbolic powers, see \cite{SurveySP}.

Various invariants have been defined to compare symbolic and ordinary powers of ideals:
the resurgence \cite{BoH}, the Waldschmidt constant  \cite{BoH}, and symbolic defect \cite{SymbolicDefect}, among others. These can be in some cases explicitly computed and in others approximated using the \texttt{SymbolicPowers} package.

\section{Basic Usage} \label{howto}

The main method in the \texttt{SymbolicPowers} package is \texttt{symbolicPower}, which given an ideal $I$ and an integer $n$ returns $I^{(n)}$. Computations are done using the standard definition of symbolic powers; if the option \texttt{UseMinimalPrimes} is set true, then the definition of symbolic powers used in the computations will be the non-standard one, as described in the introduction. When \texttt{UseMinimalPrimes} is set true, the algorithm takes a primary decomposition of $I^n$ and intersects the components corresponding to minimal primes. Through the rest of the paper, we will assume that the \texttt{UseMinimalPrimes} option is set to false, which is the default setting.

The fastest algorithm used in \texttt{symbolicPower} is the one for homogeneous ideals $I$ of height $\textrm{dim} \left( R \right) - 1$ in $R$. In this case, the $I^{(n)}$ coincides with the saturation of $I^n$ with respect to the maximal ideal, as can be seen in the following example.
\begin{example}[Height of ideal is $\textrm{dim} \left( R \right) - 1$] $\ $
\begin{verbatim}
    i1 : loadPackage "SymbolicPowers"

    i2 : R=QQ[x,y,z];

    i3 : I=ideal(x*(y^3-z^3),y*(z^3-x^3),z*(x^3-y^3));
    o3 : ideal of R

    i4 : symbolicPower(I,3);
    o4 : ideal of R

    i4 : symbolicPower(I,3)==saturate(I^3)
    o4 : True
\end{verbatim}
\end{example}

If $I$ is a primary ideal, a primary decomposition of $I^n$ is scanned for the component with radical $I$.
\begin{example}[Primary ideals] $\ $
\begin{verbatim}
    i1 : loadPackage "SymbolicPowers";

    i2 : R=QQ[w,x,y,z]/(x*y-z^2);

    i3 : I=ideal(x,z);
    o3 : Ideal of R

    i4 : symbolicPower(I,2)
    o4 = ideal x
    o4 : Ideal of R
\end{verbatim}
\end{example}

If $I$ is a monomial ideal, the symbolic powers are computed by intersecting the powers of the associated primes in the squarefree case, or by explicitly picking out the appropriate components of the primary decomposition of $I^n$ in the general case. This does not require that $I$ be of type \texttt{monomialIdeal}, but it does make use of the fact that the \texttt{primaryDecomposition} algorithm is faster for ideals of type \texttt{monomialIdeal}. 

\pagebreak

\begin{example}[Monomial ideals] $\ $
\begin{verbatim}
    i1 : loadPackage "SymbolicPowers";

    i2 : R = QQ[x,y,z];

    i3 : I = ideal(x*y,x*z,y*z)
    o3 = ideal (x*y, x*z, y*z)
    o3 : Ideal of QQ[x, y, z]

    i4 : symbolicPower(I,2)
                        2 2   2 2   2 2
    o4 = ideal (x*y*z, y z , x z , x y )
    o4 : Ideal of QQ[x, y, z]
\end{verbatim}
\end{example}

\section{Applications} \label{examples}

\subsection{Equality} \label{equality}

Symbolic powers do not, in general, coincide with the ordinary powers, even in the case of prime ideals. In fact, the question of characterizing the ideals $I$ for which $I^{(n)} = I^n$ for all $n$ is essentially open. Using \texttt{isSymbolicEqualOrdinary}, one can determine if the $n$-th symbolic and ordinary powers of a given ideal coincide, often without computing the actual symbolic power of $I$. For this, the package makes use of \texttt{bigHeight}.

The method \texttt{bigHeight} computes the largest height of an associated prime of $I$. Similarly, \texttt{minimalPart} returns the intersection of the minimal components of a given ideal. Instead of explicitly finding the associated primes of $I$ and taking their heights, the following result is used \cite{EisenbudHunekeVasconcelosPrimaryDecomposition}.

\begin{thm}[Eisenbud--Huneke--Vasconcellos, 1992]
	Given an ideal $I$ in a regular ring $R$ of height $h$, then for each $e \geqslant h$, $I$ has an associated prime of height $e$ if and only if the height of $\Ext^e \left( R/I, R\right)$ is $e$. 
\end{thm}

This is faster than computing the set of associated primes. To determine if $I^{(n)} = I^n$ for a specific value of $n$, \texttt{isSymbolicEqualOrdinary} first compares the big heights of $I^n$ and $I$: if the big heights differ, then $I^n$ must have embedded components, and \texttt{isSymbolicEqualOrdinary} returns \texttt{false}; if the big heights are both equal to the height of $I$, then $I^n$ cannot have embedded components, and \texttt{isSymbolicEqualOrdinary} returns \texttt{true}.

\subsection{The Containment Problem} \label{containment}

The Containment Problem for ordinary and symbolic powers of ideals consists of answering the following question: given an ideal $I$, for which values of $a$ and $b$ does the containment $I^{(a)} \subseteq I^b$ hold? Over a regular ring, a well-known theorem of Ein--Lazersfeld--Smith, Hochster--Huneke and Ma--Schwede \cite{ELS,comparison,MaSchwedeSymbPowers} gives a partial answer to this question: when $I$ is a radical ideal, $I^{(hn)} \subseteq I^n$ holds for all $n$, where $h$ denotes the big height of the ideal $I$. However, this is not necessarily best possible. Using \texttt{containmentProblem}, the user can either determine the smallest value of $a$ given $b$ for which $I^{(a)} \subseteq I^b$, or the largest value of $b$ given $a$ for which the same containment holds.
\begin{example}[Containment Problem] $\ $
\begin{verbatim}
    i1 : loadPackage "SymbolicPowers";

    i2 : R=QQ[x,y,z];

    i3 : I=ideal(x*(y^3-z^3),y*(z^3-x^3),z*(x^3-y^3));
    o3 : ideal of R

    i4 : containmentProblem(I,2)
    o4 = 4
\end{verbatim}
\end{example}


\section{Asymptotic invariants} \label{asymptotic}
In an effort to make progress on the containment problem, various asymptotic interpolation invariants have been  proposed by Bocci and Harbourne \cite{BoH}. One such invariant is  the Waldschmidt constant for a homogeneous ideal $I$.  This is an asymptotic measure of the initial degree of the symbolic powers of $I$. The  {\em initial degree} of a homogeneous ideal $I$ is $\alpha(I)~=~\min\{ d ~|~ I_d  \neq 0 \}$, i.e. the smallest degree of a nonzero element in $I$. The {\em Waldschmidt constant} of $I$ is defined to be
$$\widehat\alpha(I) = \lim_{m \rightarrow \infty} \frac{\alpha(I^{(m)})}{m}. $$

Due to the asymptotic nature of the Waldschmidt constant, there is no a priori algorithm to determine this for arbitrary ideals. An important exception is the case when the ideal $I$ is a monomial ideal. In this context, the  Waldschmidt constant can be computed as the smallest among the sums of the coordinates of all points in a convex body termed the {\em symbolic polyhedron} of $I$ \cite{SymbMon,WaldschmidtMonomials}. Our package computes Waldschmidt constants of monomial ideals by finding their symbolic polyhedron. The \texttt{symbolicPolyhedron} routine makes heavy use of the {\em Polyhedra} package by Ren\'e Birkner, which in turn relies on the {\em FourierMotzkin} package by Greg Smith. This allows to determine the Waldschmidt constants of monomial ideals exactly as in the following example.

\begin{example}[Waldschmidt constant of monomial ideals] $\ $
\begin{verbatim}
    i1 : loadPackage "SymbolicPowers";

    i2 : R=QQ[x,y,z];

    i3 : I=ideal(x*y,x*z,y*z);
    o3 : Ideal of R

    i4 : symbolicPolyhedron(I)
    o4 = {ambient dimension => 3           }
          dimension of lineality space => 0
          dimension of polyhedron => 3
          number of facets => 6
          number of rays => 3
          number of vertices => 4
    o4 : Polyhedron

    i5 : waldschmidt I
    Ideal is monomial, the Waldschmidt constant is computed exactly
\end{verbatim}

\pagebreak 

\begin{verbatim}
          3
    o5 =  - 
          2
    o5 : QQ
\end{verbatim}
\end{example}

In the case of arbitrary ideals, the Waldschmidt constant is approximated by taking the minimum of the values 
$\frac{\alpha(I^{(m)})}{m}$, where $m$ ranges from 1 to a specified optional input \texttt{SampleSize}. 

\begin{example}[Waldschmidt constant of arbitrary ideals] $\ $
\begin{verbatim}
    i1 : loadPackage "SymbolicPowers";

    i2 : R=QQ[x,y,z];

    i3 : I=ideal(x*(y^3-z^3),y*(z^3-x^3),z*(x^3-y^3));
    o3 : Ideal of R

    i4 : waldschmidt I
    Ideal is not monomial, the  Waldschmidt constant is approximated 
    using first 10 powers.
    o4 = 3
    o4 : QQ
\end{verbatim}
\end{example}

Note that the true value for the Waldschmidt constant of the above ideal is indeed $3$ as proven in \cite{DHNSST2015}. In general, for an ideal that is not monomial, the function \texttt{waldschmidt} will return an upper bound on the true value of the Waldschmidt constant.

Another asymptotic invariant termed {\em resurgence} is defined as
$$\rho(I) = \sup\left\{\frac{m}{r} ~|~ I^{(m)} \not\subseteq I^r \right\}.$$
There are no algorithms known to date that compute resurgence exactly, therefore our package computes a lower bound for the resurgence by taking the maximum of the values $\frac{m}{r}$, where $r$ ranges from 1 to the second input of the function \texttt{lowerBoundResurgence}. Continuing with the ideal in the previous example, we compute a lower bound on its resurgence.

\begin{example}[Lower bound on resurgence] $\ $
\begin{verbatim}
    i5 : lowerBoundResurgence(I,5)
         3
    o5 = -
         2
    o5 : QQ
\end{verbatim}
\end{example}
Note that  the value for the resurgence of the ideal in this example is also $\frac{3}{2}$ by \cite{DHNSST2015}.\section*{Acknowledgments}

We would like to thank the organizers of the July 2017 Macaulay2 Workshop at the University of California -- Berkeley, where a large portion of this work was done. The code for computing the symbolic polyhedron and Waldschmidt constant of a monomial ideal was developed by the third named author in collaboration with Andrew Conner and Xuehua (Diana) Zhong. We thank them for their contribution to these routines.

\bibliographystyle{alpha}
\bibliography{References}
\end{document}